\documentclass[preprint,12pt]{elsarticle}

\usepackage{amssymb}
\newtheorem{theorem}{Theorem}[section]
\newtheorem{lemma}[theorem]{Lemma}
\newtheorem{corollary}[theorem]{Corollary}

\newtheorem{definition}[theorem]{Definition}

\newtheorem{remark}[theorem]{Remark}

\journal{}

\begin{document}

\begin{frontmatter}



\title{Borel reducibility and finitely H\" older($\alpha$) embeddability}


\author{Longyun Ding\footnote{Research partially supported by the National Natural Science Foundation of China (Grant No. 10701044).
We thank Rui Liu for the inspiring discussions.}}

\ead{dinglongyun@gmail.com}
\address{School of Mathematical Sciences and LPMC, Nankai University, Tianjin 300071, PR China}

\begin{abstract}

Let $(X_n,d_n),\,n\in\Bbb N$ be a sequence of pseudo-metric spaces, $p\ge 1$. For $x,y\in\prod_{n\in\Bbb N}X_n$,
let $(x,y)\in E((X_n)_{n\in\Bbb N};p)\Leftrightarrow\sum_{n\in\Bbb N}d_n(x(n),y(n))^p<+\infty$. For Borel reducibility between equivalence relations
$E((X_n)_{n\in\Bbb N};p)$, we show it is closely related to finitely H\"older($\alpha$) embeddability between pseudo-metric spaces.

\end{abstract}


\end{frontmatter}



\section{Introduction}

A topological space is called a {\it Polish space} if it is homeomorphic to a separable complete metric space. Let $X,Y$ be Polish spaces and
$E,F$ equivalence relations on $X,Y$ respectively. A {\it Borel
reduction} from $E$ to $F$  is a Borel function $\theta:X\to Y$ such
that
$$(x,y)\in E\iff(\theta(x),\theta(y))\in F$$ for all $x,y\in X$. We
say that $E$ is {\it Borel reducible} to $F$, denoted $E\le_B F$, if
there is a Borel reduction from $E$ to $F$. If $E\le_B F$ and
$F\le_B E$, we say that $E$ and $F$ are {\it Borel bireducible} and
denote $E\sim_B F$. We refer to \cite{BK} and \cite{gaobook}
for background on Borel reducibility.

It was proved by R. Dougherty and G. Hjorth \cite{DH} that, for $p,q\ge 1$,
$$\Bbb R^\Bbb N/\ell_p\le_B\Bbb R^\Bbb N/\ell_q\iff p\le q.$$
The equivalence relation $\Bbb R/\ell_p$ was extended to so called $\ell_p$-like equivalence relations in \cite{ding1}.
Let $(X_n,d_n),\,n\in\Bbb N$ be a sequence of pseudo-metric spaces, $p\ge 1$. For $x,y\in\prod_{n\in\Bbb N}X_n$,
$(x,y)\in E((X_n)_{n\in\Bbb N};p)\Leftrightarrow\sum_{n\in\Bbb N}d_n(x(n),y(n))^p<+\infty$.

A special case concerning separable Banach
spaces was investigated in \cite{ding}. It was showed in \cite{ding} that Borel reducibility between this kind
of equivalence relations is related to the existence of H\"older($\alpha$) embeddings.
In this paper, we introduce the notion of $C$-finitely H\"older($\alpha$) embeddability, and generalize the connection between Borel reducibility and
finitely H\"older($\alpha$) embeddability to a rather general type of metric spaces.

\section{$\ell_p$-like equivalence relations on pseudo-metric spaces}

\begin{definition}
Let $(X_n,d_n),\,n\in\Bbb N$ be a sequence of pseudo-metric spaces, $p\ge 1$. We define an equivalence relation
$E((X_n,d_n)_{n\in\Bbb N};p)$ on $\prod_{n\in\Bbb N}X_n$ by
$$(x,y)\in E((X_n,d_n)_{n\in\Bbb N};p)\iff\sum_{n\in\Bbb
N}d_n(x(n),y(n))^p<+\infty$$ for $x,y\in\prod_{n\in\Bbb N}X_n$.
We call it an {\it $\ell_p$-like equivalence relation}.
\end{definition}
If $(X_n,d_n)=(X,d)$ for every $n\in\Bbb N$, we write
$E((X,d);p)=E((X_n,d_n)_{n\in\Bbb N};p)$ for the sake of brevity. If there is no danger of confusion, we simply write
$E((X_n)_{n\in\Bbb N};p)$ and $E(X;p)$ instead of $E((X_n,d_n)_{n\in\Bbb N};p)$ and $E((X,d);p)$.

\begin{definition}
If $X$ is a Polish space, $d$ is a Borel pseudo-metric on $X$, we say $(X,d)$ is a Borel pseudo-metric space.
\end{definition}

Let $(Y_n,\delta_n),\,n\in\Bbb N$ be a sequence of pseudo-metric spaces, $y^*\in\prod_{n\in\Bbb N}Y_n$.
For $q\ge 1$, we denote by $\ell_q((Y_n)_{n\in\Bbb N},y^*)$ the pseudo-metric space whose
underlying space is
$$\left\{y\in\prod_{n\in\Bbb N}Y_n:\sum_{n\in\Bbb N}\delta_n(y(n),y^*(n))^q<+\infty\right\},$$
with the pseudo-metric
$$\delta_q(x,y)=\left(\sum_{n\in\Bbb N}\delta_n(x(n),y(n))^q\right)^\frac{1}{q}.$$

\begin{theorem}\label{reduce}
Let $(Y,\delta)$ be a Borel pseudo-metric space, $Y_0\subseteq Y_1\subseteq Y_2\subseteq\cdots$ a sequence of Borel subsets of $Y$,
and let $(X_n,d_n),n\in\Bbb N$ be a sequence of Borel
pseudo-metric spaces, $p,q\in[1,+\infty)$. If there are $A,C,D>0$, a sequence
of Borel maps $T_n:X_n\to\ell_q((Y_n)_{n\in\Bbb N},y^*)$ for some $y^*\in\prod_{n\in\Bbb N}Y_n$
and two sequences of non-negative real numbers
$\varepsilon_n,\eta_n,n\in\Bbb N$ such that
\begin{enumerate}
\item[(1)] $\sum_{n\in\Bbb N}\varepsilon_n^p<+\infty,\sum_{n\in\Bbb N}\eta_n^q<+\infty$;
\item[(2)]
$d_n(u,v)<\varepsilon_n\Rightarrow \delta_q(T_n(u),T_n(v))<\eta_n$;
\item[(3)] $d_n(u,v)\ge C\Rightarrow \delta_q(T_n(u),T_n(v))\ge D$;
\item[(4)] $\varepsilon_n\le d_n(u,v)<C\Rightarrow A^{-1}d_n(u,v)^{\frac{p}{q}}\le \delta_q(T_n(u),T_n(v))\le
Ad_n(u,v)^{\frac{p}{q}}$.
\end{enumerate}
Then we have $$E((X_n)_{n\in\Bbb N};p)\le_B E((Y_n)_{n\in\Bbb N};q).$$
\end{theorem}

{\bf Proof.} Fix a bijection
$\langle\cdot,\cdot\rangle:\Bbb N^2\to\Bbb N$ such that $m\le\langle n,m\rangle$ for each $n,m\in\Bbb N$. Note that
$T_n(u)(m)\in Y_m\subseteq Y_{\langle n,m\rangle}$ for every $u\in X_n$. We define $\theta:\prod_{n\in\Bbb N}X_n\to \prod_{k\in\Bbb N}Y_k$ by
$$\theta(x)(\langle n,m\rangle)=T_n(x(n))(m)$$
for $x\in\prod_{n\in\Bbb N}X_n$ and $n,m\in\Bbb N$. It is easy to
see that $\theta$ is Borel. By the definition we have
$$\begin{array}{ll}&\sum_{n,m\in\Bbb N}\delta(\theta(x)(\langle n,m\rangle),\theta(y)(\langle n,m\rangle))^q\cr
=&\sum_{n\in\Bbb N}\sum_{m\in\Bbb N}\delta(T_n(x(n))(m),T_n(y(n))(m))^q\cr
=&\sum_{n\in\Bbb N}\delta_q(T_n(x(n)),T_n(y(n)))^q\, .\end{array}$$

For $x,y\in\prod_{n\in\Bbb N}X_n$, we split $\Bbb N$ into three sets
$$I_1=\{n\in\Bbb N:d_n(x(n),y(n))<\varepsilon_n\},$$
$$I_2=\{n\in\Bbb N:d_n(x(n),y(n))\ge C\},$$
$$I_3=\{n\in\Bbb N:\varepsilon_n\le d_n(x(n),y(n))<C\}.$$
From (2) we have $$\sum_{n\in I_1}d_n(x(n),y(n))^p<\sum_{n\in
I_1}\varepsilon_n^p\le\sum_{n\in\Bbb N}\varepsilon_n^p<+\infty,$$
$$\sum_{n\in I_1}\delta_q(T_n(x(n)),T_n(y(n)))^q<\sum_{n\in I_1}\eta_n^q\le\sum_{n\in\Bbb N}\eta_n^q<+\infty;$$
denote $|I_2|$ the cardinal of $I_2$, from (3) we have
$$\sum_{n\in I_2}d_n(x(n),y(n))^p\ge C^p|I_2|,$$
$$\sum_{n\in I_2}\delta_q(T_n(x(n)),T_n(y(n)))^q\ge D^q|I_2|;$$
and from (4) we have
$$A^{-q}\sum_{n\in I_3}d_n(x(n),y(n))^p\le\sum_{n\in I_3}\delta_q(T_n(x(n)),T_n(y(n)))^q\le A^q\sum_{n\in I_3}d_n(x(n),y(n))^p.$$
Therefore,
$$\begin{array}{ll}&(x,y)\in E((X_n)_{n\in\Bbb N};p)\cr
\iff &\sum_{n\in\Bbb N}d_n(x(n),y(n))^p<+\infty\cr
\iff &|I_2|<\infty,\,\sum_{n\in I_3}d_n(x(n),y(n))^p<+\infty\cr
\iff &|I_2|<\infty,\,\sum_{n\in I_3}\delta_q(T_n(x(n)),T_n(y(n)))^q<+\infty\cr
\iff &\sum_{n\in\Bbb N}\delta_q(T_n(x(n)),T_n(y(n)))^q<+\infty\cr
\iff &\sum_{n,m\in\Bbb N}\delta(\theta(x)(\langle n,m\rangle),\theta(y)(\langle n,m\rangle))^q<+\infty\cr
\iff &(\theta(x),\theta(y))\in E((Y_k)_{k\in\Bbb N};q).$$
\end{array}$$
It follows that $E((X_n)_{n\in\Bbb N};p)\le_B E((Y_n)_{n\in\Bbb N};q)$. \hfill$\Box$

\begin{corollary} \label{reduce1}
If all $(X_n,d_n)$'s are separable, then the sequence $\eta_n,n\in\Bbb N$ and clause (2) in Theorem
\ref{reduce} can be omitted.
\end{corollary}

{\bf Proof.} By Zorn's lemma, we can find a set $S_n\subseteq X_n$
for each $n$ such that
\begin{enumerate}
\item[(i)] $\forall r,s\in S_n(r\ne s\to d_n(r,s)\ge\varepsilon_n)$;
\item[(ii)] $\forall u\in X_n\exists s\in
S_n(d_n(u,s)<\varepsilon_n)$.
\end{enumerate}
Since $X_n$ is separable, $S_n$ is countable. So we can enumerate
$S_n$ by $(s^n_m)_{m\in\Bbb N}$. Define $T'_n:X_n\to\ell_q((Y_n)_{n\in\Bbb N},y^*)$ by
$T'_n(u)=T_n(s^n_{m(u)})$ where $m(u)$ is the least $m$ such that
$d_n(u,s^n_m)<\varepsilon_n$. It is easy to see that each $T'_n$ is
Borel.

Without loss of generality, we may assume that $5\varepsilon_n<C$.
Now denote
$\varepsilon_n'=3\varepsilon_n,\eta_n'=A(5\varepsilon_n)^\frac{p}{q}$
and
$A'=3^\frac{p}{q}A,C'=C-2\varepsilon_n,D'=\min\left\{D,A^{-1}\left(\frac{C}{5}\right)^\frac{p}{q}\right\}$.
We check that $\varepsilon_n',\eta_n',A',C'$ and $D'$ meet clauses (1)--(4) in Theorem \ref{reduce} as follows:

(1) $\sum_{n\in\Bbb N}(\varepsilon_n')^p=3^p\sum_{n\in\Bbb N}\varepsilon_n^p<+\infty$, $\sum_{n\in\Bbb N}(\eta_n')^q=5^pA^q\sum_{n\in\Bbb N}\varepsilon_n^p<+\infty$.

(2) If $d_n(u,v)<\varepsilon_n'$, then
$d_n(s^n_{m(u)},s^n_{m(v)})<5\varepsilon_n<C$. Note that $s^n_{m(u)}=s^n_{m(v)}$ or $d_n(s^n_{m(u)},s^n_{m(v)})\ge\varepsilon_n$.
So by clause (4) in Theorem \ref{reduce}, we have
$$\delta_q(T_n'(u),T_n'(v))=\delta_q(T_n(s^n_{m(u)}),T_n(s^n_{m(v)}))\le
Ad_n(s^n_{m(u)},s^n_{m(v)})^\frac{p}{q}<\eta_n'.$$

(3) If $d_n(u,v)\ge C'$, then $d_n(s^n_{m(u)},s^n_{m(v)})\ge
C-4\varepsilon_n\ge\varepsilon_n$. For $\varepsilon_n\le
d_n(s^n_{m(u)},s^n_{m(v)})<C$, we have
$$\begin{array}{ll}\delta_q(T_n'(u),T_n'(v))&=\delta_q(T_n(s^n_{m(u)}),T_n(s^n_{m(v)}))\cr
&\ge A^{-1}d_n(s_{m(u)}^n,s_{m(v)}^n)^\frac{p}{q}\ge
A^{-1}(C-4\varepsilon_n)^\frac{p}{q}\cr &\ge
A^{-1}\left(\frac{C}{5}\right)^\frac{p}{q}\ge D'.\end{array}$$ And
for $d_n(s^n_{m(u)},s^n_{m(v)})\ge C$, we have
$$\delta_q(T_n'(u),T_n'(v))=\delta_q(T_n(s^n_{m(u)}),T_n(s^n_{m(v)}))\ge D\ge D'.$$

(4) If $\varepsilon_n'\le d_n(u,v)<C'$, then $\varepsilon_n\le
d_n(s^n_{m(u)},s^n_{m(v)})<C$ and
$$\frac{1}{3}d_n(u,v)\le d_n(u,v)-2\varepsilon_n<d_n(s^n_{m(u)},s^n_{m(v)})<d_n(u,v)+2\varepsilon_n\le 3d_n(u,v).$$
Since $$A^{-1}d_n(s^n_{m(u)},s^n_{m(v)})^{\frac{p}{q}}\le
\delta_q(T_n(s^n_{m(u)}),T_n(s^n_{m(v)}))\le
Ad_n(s^n_{m(u)},s^n_{m(v)})^{\frac{p}{q}},$$ it follows that
$$(A')^{-1}d_n(u,v)^\frac{p}{q}\le \delta_q(T_n'(u),T_n'(v))\le A'd_n(u,v)^\frac{p}{q}.$$
\hfill$\Box$

\section{On separable pseudo-metric spaces}

For the rest of this paper, we focus on such $E((X_n,d_n)_{n\in\Bbb N};p)$ that all $(X_n,d_n)$'s are separable Borel pseudo-metric spaces.

Let $S_n=\{s^n_m:m\in\Bbb N\}$ be a countable dense subset of $X_n$. We may assume that $d_n(s^n_{m},s^n_{k})>0$ for $m\ne k$, i.e. $(S_n,d_n)$ is a countable metric
space. For $u\in X_n$, let $m_n(u)=\min\{m:d_n(u,s^n_m)<2^{-n}\}$ and
$\vartheta:\prod_{n\in\Bbb N}X_n\to\prod_{n\in\Bbb N}D_n$ as $\vartheta(x)(n)=s^n_{m_n(u)}$.
Since $\sum_{n\in\Bbb N}d_n(x(n),\vartheta(x)(n))^p<\sum_{n\in\Bbb N}2^{-np}<+\infty$, we have $(x,\vartheta(x))\in E((X_n)_{n\in\Bbb N};p)$.
Thus $\vartheta$ is a Borel reduction of $E((X_n)_{n\in\Bbb N};p)$ to $E((S_n)_{n\in\Bbb N};p)$. So $E((X_n)_{n\in\Bbb N};p)\sim_B E((S_n)_{n\in\Bbb N};p)$.
Now let $(\overline{S_n},\overline{d_n})$ be the completion of $(S_n,d_n)$. Since $(\overline{S_n},\overline{d_n})$ is a Polish space, by the same arguments, we have
$$E((\overline{S_n})_{n\in\Bbb N};p)\sim_B E((S_n)_{n\in\Bbb N};p)\sim_B E((X_n)_{n\in\Bbb N};p).$$

Therefore, from now on, we may assume that all $(X_n,d_n)$'s are separable complete metric space.

\begin{definition}
Let $(X,d)$ be a separable complete metric space, $(F_n)_{n\in\Bbb N}$ a sequence of finite subsets of $X$. If
$F_0\subseteq F_1\subseteq\cdots\subseteq F_n\subseteq\cdots$ and $\bigcup_{n\in\Bbb N}F_n$ is dense in $X$, then we denote
$$F(X;p)=E((F_n)_{n\in\Bbb N};p).$$
\end{definition}

The following lemma shows that, under Borel bireducibility, $F(X;p)$ is independent to the choice of $(F_n)_{n\in\Bbb N}$.

\begin{lemma}
Let $(X,d)$ be a separable complete metric space, and let $(F_n)_{n\in\Bbb N}$ and $(F_n')_{n\in\Bbb N}$ be two sequences of finite subsets of $X$ satisfying that
$$F_0\subseteq F_1\subseteq\cdots\subseteq F_n\subseteq\cdots,\quad F_0'\subseteq F_1'\subseteq\cdots\subseteq F_n'\subseteq\cdots,$$
and both $\bigcup_{n\in\Bbb N}F_n$ and $\bigcup_{n\in\Bbb N}F_n'$ are dense in $X$. Then for each $p\ge 1$we have
$$E((F_n)_{n\in\Bbb N};p)\sim_B E((F_n')_{n\in\Bbb N};p).$$
\end{lemma}

{\bf Proof.} It will suffice to show that $E((F_n)_{n\in\Bbb N};p)\le_B E((F_n')_{n\in\Bbb N};p)$.
For $k\in\Bbb N$, let $\gamma_k=\min\{d(u,v):u,v\in F_k,u\ne v\}$.
Note that $\bigcup_{n\in\Bbb N}F_n'$ is dense in
$X$. For $u\in F_k$, we can find a $T_k(u)\in\bigcup_{n\in\Bbb N}F_n'$ such that $d(u,T_k(u))<\gamma_k/4$.
Then for distinct $u,v\in F_k$ we have
$$\frac{1}{2}d(u,v)\le d(u,v)-\gamma_k/2<d(T_k(u),T_k(v))<d(u,v)+\gamma_k/2\le 2d(u,v).$$
Since $F_k$ is finite, there is $n_k$ such that
$T_k(u)\in F_{n_k}'$ for each $u\in F_k$. We may assume that $(n_k)_{k\in\Bbb N}$ is strictly increasing. Fix a point
$u_0\in F_0'\subseteq F_n'$. We define $\theta:\prod_{n\in\Bbb N}F_n\to\prod_{n\in\Bbb N}F_n'$ by
$$\theta(x)(n)=\left\{\begin{array}{ll}T_k(x(k)), & n=n_k\cr u_0, &\mbox{otherwise}.\end{array}\right.$$
Then for $x,y\in\prod_{n\in\Bbb N}F_n$ we have
$$\frac{1}{2^p}\sum_{k\in\Bbb N}d(x(k),y(k))^p\le\sum_{n\in\Bbb N}d(\theta(x)(n),\theta(y)(n))^p\le 2^p\sum_{k\in\Bbb N}d(x(k),y(k))^p,$$
It follows that $\theta$ is a Borel reduction of $E((F_n)_{n\in\Bbb N};p)$ to $E((F_n')_{n\in\Bbb N};p)$.
\hfill$\Box$

\begin{remark}
We can see that $E(X;p)\sim_B F(X;p)$ when $X$ is compact. But whether it is always true for every separable complete metric space? We do not know the answer.
\end{remark}

\begin{definition}
For two metric spaces $(X,d),(X',d')$ and $\alpha>0$. We say that $X$ {\it H\"older$(\alpha)$ embeds} into $X'$ if there
exist $A>0$ and $T:X\to X'$ such that, for $u,v\in F$,
$$A^{-1}d(u,v)^\alpha\le d'(T(u),T(v))\le Ad(u,v)^\alpha.$$
\end{definition}

Theorem \ref{reduce} gives the following result.

\begin{remark}\label{reduce2}
Let $X,Y$ be two separable complete metric spaces,
$p,q\in[1,+\infty)$. If $X$ H\"older$(\frac{p}{q})$ embeds into
$\ell_q(Y,y^*)$ for some $y^*\in Y^\Bbb N$, then we have
$E(X;p)\le_B E(Y;q)$.
\end{remark}

In next section, we present a necessary condition of $E(X;p)\le_B E(Y;q)$ which will be named finitely
H\"older$(\frac{p}{q})$ embeddability.

\section{Finitely H\" older$(\alpha)$ embeddability}

A weak version of the following lemma is due to R. Dougherty and G. Hjorth \cite{DH}. For self-contain reason, we present a proof for it.

\begin{lemma}\label{DH}
Let $(Y_n,\delta_n),n\in\Bbb N$ be a a sequence of separable complete metric space, $p,q\in[1,+\infty)$,
and let $(Z_n,d_n),n\in\Bbb N$ be a sequence of finite metric spaces.
Assume that $E((Z_n)_{n\in\Bbb N};p)\le_B E((Y_n)_{n\in\Bbb N};q)$. Then there exist
strictly increasing sequences of natural numbers $(b_j)_{j\in\Bbb
N},(l_j)_{j\in\Bbb N}$ and $T_j:Z_{b_j}\to\prod_{n=l_j}^{l_{j+1}-1}Y_n$ such
that, for $x,y\in\prod_{j\in\Bbb N}Z_{b_j}$, we have
$$(x,y)\in E((Z_{b_j},d_{b_j})_{j\in\Bbb N};p)\iff\sum_{j\in\Bbb N}\delta_q(T_j(x(j)),T_j(y(j)))^q<+\infty,$$
where $\delta_q(r,s)=(\sum_{n=l_j}^{l_{j+1}-1}\delta_n(r(n),s(n))^q)^\frac{1}{q}$ for $r,s\in \prod_{n=l_j}^{l_{j+1}-1}Y_n$.
\end{lemma}

{\bf Proof.} The proof is modified from the proof of \cite{DH} Theorem 2.2, Claim (i)--(iii).

Denote $Z=\prod_{n\in\Bbb N}Z_n$. Assume that $\theta$ is a Borel reduction of $E((Z_n)_{n\in\Bbb N};p)$ to $E((Y_n)_{n\in\Bbb N};q)$.
For each finite sequence $t$ we denote $l(t)$ the length of $t$; if $t\in\prod_{i<l(t)}Z_i$, let $N_t=\{z\in Z:z(i)=t(i)\,(i<l(t))\}$.

{\it Claim} (i). For $j,k\in\Bbb N$, there exist $l\in\Bbb N$ and $s^*\in\prod_{i=k}^{k+l(s^*)-1}Z_i$ and a comeager set $D\subseteq Z$ such that,
for all $x,\hat x\in D$, if we have $x=rs^*y$ and $\hat x=\hat rs^*y$ for some $r,\hat r\in\prod_{i<k}Z_i$ and $y\in\prod_{i\ge k+l(s^*)}Z_i$, then
$$\sum_{n\ge l}\delta_n(\theta(x)(n),\theta(\hat x)(n))^q<2^{-j}.$$

{\it Proof.} For $l\in\Bbb N$, we define a function $F_l:Z\to\Bbb R$ by
$$F_l(x)=\max\left\{\sum_{n\ge l}\delta_n(\theta(z)(n),\theta(\hat z)(n))^q:z(i)=\hat z(i)=x(i)\,(i\ge k)\right\}.$$
For each $x$, there are only finitely many pairs $z,\hat z$ satisfying $z(i)=\hat z(i)=x(i)\,(i\ge k)$. For each such pair we have
$(z,\hat z)\in E((Z_n)_{n\in\Bbb N};p)$, so $(\theta(z),\theta(\hat z))\in E((Y_n)_{n\in\Bbb N};q)$. Thus
$\lim_{l\to\infty}\sum_{n\ge l}\delta_n(\theta(z)(n),\theta(\hat z)(n))^q=0$.
Hence $F_l(x)<+\infty$ for all $l$ and $\lim_{l\to\infty}F_l(x)=0$. Therefore, by the Baire category theorem, there exists an $l$ such that $\{x:F_l(x)<2^{-j}\}$ is not
meager. By $F$ is Borel, this set has the property of Baire, so there is an open set $O\ne\emptyset$ on which it is relatively comeager.

Find an $N_t\subseteq O$ for some finite sequence $t$ with $l(t)\ge k$. Let $t=r^*s^*$ where $l(r^*)=k$. Since $F_l(x)$ does not depend on the first $k$ coordinates of $x$,
we have $\{x:F_l(x)<2^{-j}\}$ is also relatively comeager in $N_{rs^*}$ for all $r\in\prod_{i<k}Z_i$. Let $D$ be a comeager set such that $F_l(x)<2^{-j}$ whenever
$x\in D\cap N_{rs^*}$ for any $r$ of length $k$. Now the conclusion of the claim follows from the definition of $F_l$.
\hfill Claim (i) $\Box$

By \cite{kechris} Theorem (5.38), there is a dense $G_\delta$ set $C\subseteq Z$ such that $\theta\upharpoonright C$ is continuous.

{\it Claim} (ii). For $j,k,l\in\Bbb N$, there exists a finite sequence $s^{**}\in\prod_{i=k}^{k+l(s^{**})-1}Z_i$ such that,
for all $x,\hat x\in C$, if we have $x=rs^{**}y$ and $\hat x=rs^{**}\hat y$ for some $r\in\prod_{i<k}Z_i$ and $y,\hat y\in\prod_{i\ge k+l(s^{**})}Z_i$, then
$$\sum_{n<l}\delta_n(\theta(x)(n),\theta(\hat x)(n))^q<2^{-j}.$$

Furthermore, if $G$ is a given dense open subset of $Z$, then $s^{**}$ can be chosen such that $N_{rs^{**}}\subseteq G$ for all $r\in\prod_{i<k}Z_i$.

{\it Proof.} Since $\prod_{i<k}Z_i$ is a finite set, we may enumerate its elements as $r_0,r_1,\cdots,r_{M-1}$. We construct
finite sequences $t_0,t_1,\cdots,t_M$ as follows.

Let $t_0=\emptyset$. Suppose that $m<M$ and we have constructed a finite sequence $t_m\in\prod_{i=k}^{k+l(t_m)-1}Z_i$. The basic open set $N_{r_mt_m}$ must meet the comeager
set $C$, so we can pick a $w\in C\cap N_{r_mt_m}$. Since $\theta$ is continuous on $C$ and $\delta_n$ is continuous on $Y_n^2$, we can find a neighborhood $O$ of $w$
such that, for all $x,\hat x\in C\cap O$,
$\sum_{n<l}\delta_n(\theta(x)(n),\theta(\hat x)(n))^q<2^{-j}$. Find an $N_{r_mt_m'}\subseteq N_{r_mt_m}\cap O$, then $t_m\subseteq t_m'$. Since $G$ is open dense, we can
further extend $t_m'$ to get $t_{m+1}$ such that $N_{r_mt_{m+1}}\subseteq G$. Once the sequences $t_m\,(m\le M)$ are constructed, $s^{**}=t_M$ fulfills the requirements.
\hfill Claim (ii) $\Box$

We now repeatedly apply Claims (i) and (ii) to define natural numbers $b_0<b_1<b_2<\cdots$ and $l_0<l_1<l_2<\cdots$, finite sequences $(s_j)_{j\in\Bbb N}$ and dense open
sets $D_i^j\subseteq Z\,(i,j\in\Bbb N)$ as follows.

Let $b_0=l_0=0$. Suppose we have constructed $b_j,l_j,D_i^{j'}\,(j'<j)$. Applying Claim (i) for this $j$ with $k=b_j+1$, we get $l_{j+1}$, a finite sequence $s_j^*$
and a comeager set $D^j$ satisfying the conclusion of Claim (i). Let $D^j_0\supseteq D^j_1\supseteq D^j_2\supseteq\cdots$ be dense open sets of $Z$ such that
$\bigcap_{i\in\Bbb N}D^j_i\subseteq D^j\cap C$. Now apply Claim (ii) for $j$ with $k=b_j+1+l(s_j^*),l=l_{j+1}$ and $G=\bigcap_{j'<j}D_j^{j'}$ to get $s_j^{**}$.
We set $s_j=s_j^*s_j^{**}$ and $b_{j+1}=b_j+l(s_j)+1$.

Denote $Z'=\prod_{j\in\Bbb N}Z_{b_j}$ and define $h:Z'\to Z$ by
$$h(x)=\langle x(0)\rangle s_0\langle x(1)\rangle s_1\langle x(2)\rangle s_2\cdots.$$
Since $s_j=s_j^*s_j^{**}$, $h(x)$ has the form $rs_j^*y$ where $l(r)=b_j+1$, and also has the form $rs_j^{**}y$ where $l(r)=b_j+l(s^*)+1$. Therefore, Claim (ii) for
$s_j^{**}$ gives $h(x)\in G=\bigcap_{j'<j}D_j^{j'}$. Hence, for any $j$, we have $h(x)\in D_i^j$ for $i>j$, so $h(x)\in D^j\cap C$. Therefore, Claims (i) and (ii)
imply that, for any $x,\hat x\in Z'$:
\begin{enumerate}
\item[(1)] if $x(b_i)=\hat x(b_i)\,(i>j)$, then $\sum_{n\ge l_{j+1}}\delta_n(\theta(h(x))(n),\theta(h(\hat x))(n))^q<2^{-j}$;
\item[(2)] if $x(b_i)=\hat x(b_i)\,(i\le j)$, then $\sum_{n<l_{j+1}}\delta_n(\theta(h(x))(n),\theta(h(\hat x))(n))^q<2^{-j}$.
\end{enumerate}

Fix a point $u_0\in Z_0\subseteq Z_{b_i}$. For $j\in\Bbb N$ we define $T_j:Z_{b_j}\to\prod_{n=l_j}^{l_{j+1}-1}Y_n$ by
$$T_j(w)=\theta(h(\langle u_0,\cdots,u_0,w,u_0,u_0,\cdots\rangle))\upharpoonright[l_j,l_{j+1})$$
with $j$ $u_0$'s before $v$. Let $\theta':Z\to\prod_{n\in\Bbb N}Y_n$,
$$\theta'(x)=T_0(x(0))T_1(x(1))T_2(x(2))\cdots.$$
Next claim shows that $\theta'$ is a Borel reduction of $E((Z_{b_j},d_{b_j})_{j\in\Bbb N};p)$ to $E((Y_n)_{n\in\Bbb N};q)$.

{\it Claim} (iii). For all $x,\hat x\in\prod_{j\in\Bbb N}Z_{b_j}$, we have
$$(x,\hat x)\in E((Z_{b_j},d_{b_j})_{j\in\Bbb N};p)\iff(\theta'(x),\theta'(\hat x))\in E((Y_n)_{n\in\Bbb N};q).$$

{\it Proof.} Note that
$$\begin{array}{ll}(x,\hat x)\in E((Z_{b_j},d_{b_j})_{j\in\Bbb N};p)&\iff(h(x),h(\hat x))\in E((Z_n,d_n)_{n\in\Bbb N};p)\cr
&\iff(\theta(h(x)),\theta(h(\hat x)))\in E((Y_n)_{n\in\Bbb N};q).\end{array}$$
It will suffice to show that $(\theta(h(x)),\theta'(x))\in E((Y_n)_{n\in\Bbb N};q)$ for any $x\in Z'$.

For any $x\in Z'$ and $j\in\Bbb N$, define $e_j(x),e_j'(x)\in Z'$ by
$$e_j(x)(i)=\left\{\begin{array}{ll}x(i), & i=j\cr u_0, & i\ne j;\end{array}\right.\quad e_j'(x)(i)=\left\{\begin{array}{ll}x(i), & i\le j\cr u_0, & i>j.\end{array}\right.$$
By (1) for $j-1$ and (2), we have
$$\sum_{n\ge l_j}\delta_n(\theta(h(e_j(x)))(n),\theta(h(e_j'(x)))(n))^q<2^{-(j-1)},$$
$$\sum_{n<l_{j+1}}\delta_n(\theta(h(x))(n),\theta(h(e_j'(x)))(n))^q<2^{-j}.$$
Thus we have
$$\begin{array}{ll}&\sum_{n=l_j}^{l_{j+1}-1}\delta_n(\theta(h(x))(n),\theta(h(e_j(x)))(n))^q\cr
\le &\sum_{n=l_j}^{l_{j+1}-1}[\delta_n(\theta(h(x))(n),\theta(h(e_j'(x)))(n))+\delta_n(\theta(h(e_j(x)))(n),\theta(h(e_j'(x)))(n))]^q\cr
\le & 2^{q-1}\left[\sum_{n=l_j}^{l_{j+1}-1}\delta_n(\theta(h(x))(n),\theta(h(e_j'(x)))(n))^q\right.\cr
&\hskip 1cm\left.+\sum_{n=l_j}^{l_{j+1}-1}\delta_n(\theta(h(e_j(x)))(n),\theta(h(e_j'(x)))(n))^q\right]\cr
\le & 2^{q-1}\left[\sum_{n<l_{j+1}}\delta_n(\theta(h(x))(n),\theta(h(e_j'(x)))(n))^q\right.\cr
&\hskip 1cm\left.+\sum_{n\ge l_j}\delta_n(\theta(h(e_j(x)))(n),\theta(h(e_j'(x)))(n))^q\right]\cr
<&2^{q-1}\cdot 3\cdot 2^{-j}.\end{array}$$
We can see that $\theta'(x)\upharpoonright[l_j,l_{j+1})=T_j(x(j))=\theta(h(e_j(x)))\upharpoonright[l_j,l_{j+1})$ for each $j\in\Bbb N$. Therefore,
$$\begin{array}{ll}&\sum_{n\in\Bbb N}\delta_n(\theta(h(x))(n),\theta'(x)(n))^q\cr
=&\sum_{j\in\Bbb N}\sum_{n=l_j}^{l_{j+1}-1}\delta_n(\theta(h(x))(n),\theta'(x)(n))^q\cr
=&\sum_{j\in\Bbb N}\sum_{n=l_j}^{l_{j+1}-1}\delta_n(\theta(h(x))(n),\theta(h(e_j(x)))(n))^q\cr
<&\sum_{j\in\Bbb N}2^{q-1}\cdot 3\cdot 2^{-j}<+\infty,\end{array}$$
as desired.
\hfill Claim (iii) $\Box$

Note that
$$\begin{array}{ll}(\theta'(x),\theta'(\hat x))\in E((Y_n)_{n\in\Bbb N};q)&\iff\sum_{j\in\Bbb N}\sum_{n=l_j}^{l_{j+1}-1}\delta_n(\theta'(x)(n),\theta'(x)(n))^q<+\infty\cr
&\iff\sum_{j\in\Bbb N}\delta_q(T_j(x(j)),T_j(y(j)))^q<+\infty.\end{array}$$
This completes the proof.
\hfill$\Box$

Let $(X,d)$ be a metric space and $C>0$. We consider the following
condition:

$({\rm link}(C))$ {\it For $\varepsilon>0$, there exists $N\ge 1$ such that,
for any $u,v\in X$ with $d(u,v)<C$, we can find $r_i\in
X,\,i=0,1,\cdots,N$ with $r_0=u,r_N=v$ and
$d(r_{i-1},r_i)<\varepsilon$ for each $i\ge 1$.}

Let $(X,d)$ and $(Y_n,\delta_n),\,n\in\Bbb N$ be separable complete metric spaces,
$p,q\in[1,+\infty)$. Assume that
\begin{enumerate}
\item[(A1)] $X$ satisfies $({\rm link}(C))$ for some $C>0$; and
\item[(A2)] $F(X;p)\le_B E((Y_n)_{n\in\Bbb N};q)$.
\end{enumerate}

Fix a sequence of finite subsets $F_n\subseteq X,n\in\Bbb N$ such
that
$$F_0\subseteq F_1\subseteq\cdots\subseteq F_n\subseteq\cdots$$
and $\bigcup_{n\in\Bbb N}F_n$ is dense in $X$.

Since $({\rm link}(C))$ holds, for $l\in\Bbb N$, there exists $N(l)\ge 1$
such that, for any $u,v\in X$ with $d(u,v)<C$, we can find
$r^l_i(u,v)\in X,\,i=0,1,\cdots,N(l)$ with
$r^l_0(u,v)=u,r^l_{N(l)}(u,v)=v$ and
$d(r^l_{i-1}(u,v),r^l_i(u,v))<2^{-l}$ for $i=1,\cdots,N(l)$. We
denote
$$Z_n=\{r^l_i(u,v):u,v\in F_n,d(u,v)<C,l\le n,i=0,1,\cdots,N(l)\}.$$

Note that $E((Z_n);p)\sim_B F(X;p)\le_B E((Y_n)_{n\in\Bbb N};q)$. Since $Z_n\subseteq X$ is a sequence of finite metric spaces, we can
find $(b_j)_{j\in\Bbb N},(l_j)_{j\in\Bbb N}$ and $T_j:Z_{b_j}\to\prod_{n=l_j}^{l_{j+1}-1}Y_n$ as in Lemma~\ref{DH}. Then we have the following
lemmas.

\begin{lemma} \label{D}
For any $C'>0$, there exists a $D>0$ such that, for sufficiently
large $j$ and $u,v\in F_{b_j}$, if $d(u,v)\ge C'$, then
$\delta_q(T_j(u),T_j(v))\ge D$.
\end{lemma}

{\bf Proof.} Assume for contradiction that, there exists a strictly
increasing sequence of natural numbers $(j_k)_{k\in\Bbb N}$ such
that there are $u_k,v_k\in F_{b_{j_k}}$ with $d(u_k,v_k)\ge C'$ and
$\delta_q(T_{j_k}(u_k),T_{j_k}(v_k))<2^{-k}$.

Now we select $x,y\in\prod_{j\in\Bbb N}Z_{b_j}$ such that
$$\left\{\begin{array}{ll}x(j)=u_k,y(j)=v_k,& j=j_k,\cr
x(j)=y(j),& \mbox{otherwise.}\end{array}\right.$$

Then we have
$$\sum_{j\in\Bbb N}d(x(j),y(j))^p=\sum_{k\in\Bbb N}d(u_k,v_k)^p\ge\sum_{k\in\Bbb N}(C')^p=+\infty,$$
so $(x,y)\notin E((Z_{b_j})_{j\in\Bbb N};p)$. On the other hand, we
have
$$\begin{array}{ll}\sum_{j\in\Bbb N}\delta_q(T_j(x(j)),T_j(y(j)))^q&=\sum_{k\in\Bbb N}\delta_q(T_{j_k}(u_k),T_{j_k}(v_k))^q\cr
&<\sum_{k\in\Bbb N}2^{-kq}\cr &<+\infty,\end{array}$$ contradicting
Lemma~\ref{DH}! \hfill$\Box$

\begin{lemma} \label{fholder}
There exists an $m\in\Bbb N$ such that $\forall k\exists N\forall
j>N$, for $u,v\in F_{b_j}$, if $k^{-1}\le d(u,v)<C$, then we have
$$2^{-m}d(u,v)^{\frac{p}{q}}\le \delta_q(T_j(u),T_j(v))\le
2^m d(u,v)^{\frac{p}{q}}.$$
\end{lemma}

{\bf Proof.} Assume for contradiction that, for every $m$, $\exists
k_m\exists^{\infty}j\exists u_j,v_j\in F_{b_j}$ such that
$k_m^{-1}\le d(u_j,v_j)<C$ but either
$$2^{-m}d(u_j,v_j)^{\frac{p}{q}}>\delta_q(T_j(u_j),T_j(v_j))$$ or
$$\delta_q(T_j(u_j),T_j(v_j))>2^m d(u_j,v_j)^{\frac{p}{q}}.$$

We define two subsets $I_1,I_2\subseteq\Bbb N$. For $m\in\Bbb N$, we
put $m\in I_1$ iff $\exists k_m\exists^{\infty}j\exists u_j,v_j\in
F_{b_j}$ satisfying that $k_m^{-1}\le d(u_j,v_j)<C$ and
$$2^{-m}d(u_j,v_j)^{\frac{p}{q}}>\delta_q(T_j(u_j),T_j(v_j));$$
and $m\in I_2$ iff $\exists k_m\exists^{\infty}j\exists u_j,v_j\in
F_{b_j}$ satisfying that $k_m^{-1}\le d(u_j,v_j)<C$ and
$$\delta_q(T_j(u_j),T_j(v_j))>2^m d(u_j,v_j)^{\frac{p}{q}}.$$

From the assumption, we can see that $I_1\cup I_2=\Bbb N$. Now we
consider the following two cases.

{\sl Case 1.} $|I_1|=\infty$. Select a finite set $J^m\subseteq\Bbb
N$ for every $m\in I_1$ and $u_j,v_j\in F_{b_j}$ for $j\in J^m$
satisfying that
\begin{enumerate}
\item[(i)] for $j\in J^m$, we have $2^{-m}d(u_j,v_j)^{\frac{p}{q}}>\delta_q(T_j(u_j),T_j(v_j))$;
\item[(ii)] $C^p\le\sum_{j\in J^m}d(u_j,v_j)^p<2C^p$;
\item[(iii)] if $m_1<m_2$, then $\max J^{m_1}<\min J^{m_2}$.
\end{enumerate}

Now we select $x,y\in\prod_{j\in\Bbb N}Z_{b_j}$ such that
$$\left\{\begin{array}{ll}x(j)=u_j,y(j)=v_j,& j\in J^m,m\in I_1,\cr
x(j)=y(j),& \mbox{otherwise.}\end{array}\right.$$ Then we have
$$\sum_{j\in\Bbb N}d(x(j),y(j))^p=\sum_{m\in I_1}\sum_{j\in
J^m}d(u_j,v_j)^p\ge\sum_{m\in I_1}C^p=+\infty,$$ so $(x,y)\notin
E((Z_{b_j})_{j\in\Bbb N};p)$. On the other hand, we have
$$\begin{array}{ll}\sum_{j\in\Bbb N}\delta_q(T_j(x(j)),T_j(y(j)))^q&=\sum_{m\in I_1}\sum_{j\in
J^m}\delta_q(T_j(u_j),T_j(v_j))^q\cr &<\sum_{m\in I_1}\sum_{j\in J^m}
2^{-mq}d(u_j,v_j)^p\cr &<2C^p\sum_{m\in I_1}\left(2^{-q}\right)^m\cr
&<+\infty,\end{array}$$ contradicting Lemma~\ref{DH}!

{\sl Case 2.} $|I_2|=\infty$. We can find a strictly increasing
sequence of natural numbers $m_l\in I_2,\,l\in\Bbb N$ such that
$m_l\ge\frac{pl}{2q}$ and $2^{m_l}\ge N(l)$ for each $l$.

We define two subsets $L_1,L_2\subseteq\Bbb N$. For $l\in\Bbb N$,
we put $l\in L_1$ iff $\exists^{\infty}j\exists u_j,v_j\in F_{b_j}$
satisfying that $k_{m_l}^{-1}\le d(u_j,v_j)<(\sqrt 2)^{-l}$ and
$$\delta_q(T_j(u_j),T_j(v_j))>2^{m_l}d(u_j,v_j)^{\frac{p}{q}};$$
and $l\in L_2$ iff $\exists^{\infty}j\exists u_j,v_j\in F_{b_j}$
satisfying that $(\sqrt 2)^{-l}\le d(u_j,v_j)<C$ and
$$\delta_q(T_j(u_j),T_j(v_j))>2^{m_l}d(u_j,v_j)^{\frac{p}{q}}.$$

Since each $m_l\in I_2$, we have $L_1\cup L_2=\Bbb N$. We consider
two subcases.

{\sl Subcase 2.1.} $|L_1|=\infty$. Select a finite set
$K_1^l\subseteq\Bbb N$ for every $l\in L_1$ and $u_j,v_j\in F_{b_j}$
for $j\in K_1^l$ satisfying that
\begin{enumerate}
\item[(i)] for $j\in K_1^l$, we have $\delta_q(T_j(u_j),T_j(v_j))>2^{m_l}
d(u_j,v_j)^{\frac{p}{q}}$;
\item[(ii)] $(\sqrt 2)^{-pl}\le\sum_{j\in K_1^l}d(u_j,v_j)^p<2(\sqrt 2)^{-pl}$;
\item[(iii)] if $l_1<l_2$, then $\max K_1^{l_1}<\min K_1^{l_2}$.
\end{enumerate}

Now we select $x,y\in\prod_{j\in\Bbb N}Z_{b_j}$ such that
$$\left\{\begin{array}{ll}x(j)=u_j,y(j)=v_j,& j\in K_1^l,l\in L_1,\cr
x(j)=y(j),& \mbox{otherwise.}\end{array}\right.$$ Then we have
$$\sum_{j\in\Bbb N}d(x(j),y(j))^p=\sum_{l\in L_1}\sum_{j\in
K_1^l}d(u_j,v_j)^p\le 2\sum_{l\in L_1}(\sqrt 2)^{-pl}<+\infty,$$ so
$(x,y)\in E((Z_{b_j})_{j\in\Bbb N};p)$. On the other hand, since
$m_l\ge\frac{pl}{2q}$, we have
$$\begin{array}{ll}\sum_{j\in\Bbb N}\delta_q(T_j(x(j)),T_j(y(j)))^q&=\sum_{l\in L_1}\sum_{j\in
K_1^l}\delta_q(T_j(u_j),T_j(v_j))^q\cr &>\sum_{l\in L_1}\sum_{j\in K_1^l}
2^{qm_l}d(u_j,v_j)^p\cr &\ge \sum_{l\in L_1}(\sqrt 2)^{2qm_l-pl}\cr
&=+\infty,\end{array}$$ contradicting Lemma~\ref{DH}!

{\sl Subcase 2.2} $|L_2|=\infty$. Select a finite set
$K_2^l\subseteq\Bbb N$ for each $l\in L_2$ and $u_j,v_j\in F_{b_j}$
for $j\in K_2^l$ satisfying that
\begin{enumerate}
\item[(i)] for $j\in K_2^l$, we have $(\sqrt 2)^{-l}\le d(u_j,v_j)<C$ and $\delta_q(T_j(u_j),T_j(v_j))>2^{m_l}
d(u_j,v_j)^{\frac{p}{q}}$;
\item[(ii)] $C^p\le\sum_{j\in K_2^l}d(u_j,v_j)^p<2C^p$;
\item[(iii)] if $l_1<l_2$, then $\max K_2^{l_1}<\min K_2^{l_2}$;
\item[(iv)] for $j\in K_2^l$, we have $l\le b_j$.
\end{enumerate}

For $l\in L_1$ and $j\in K_2^l$, since $d(u_j,v_j)<C$ and $l\le
b_j$, by the definition of $Z_{b_j}$ we have
$$r^l_i(u_j,v_j)\in Z_{b_j}\quad(i=0,1,\cdots,N(l)).$$
Since $r^l_0(u_j,v_j)=u_j,r^l_{N(l)}(u_j,v_j)=v_j$, the triangle
inequality gives
$$\sum_{1\le i\le
N(l)}\delta_q(T_j(r^l_{i-1}(u_j,v_j)),T_j(r^l_i(u_j,v_j)))\ge
\delta_q(T_j(u_j),T_j(v_j)),$$ thus there is an $i(j)$ such that
$$\delta_q(T_j(r^l_{i(j)-1}(u_j,v_j)),T_j(r^l_{i(j)}(u_j,v_j)))\ge N(l)^{-1}\delta_q(T_j(u_j),T_j(v_j)).$$

Now denote $r_j=r^l_{i(j)-1}(u_j,v_j),s_j=r^l_{i(j)}(u_j,v_j)$. We
select $x,y\in\prod_{j\in\Bbb N}Z_{b_j}$ such that
$$\left\{\begin{array}{ll}x(j)=r_j,y(j)=s_j,& j\in K_2^l,l\ge 1,\cr
x(j)=y(j),& \mbox{otherwise.}\end{array}\right.$$ Note that
$d(r_j,s_j)<2^{-l}\le(\sqrt 2)^{-l}d(u_j,v_j)$, we have
$$\begin{array}{ll}\sum_{j\in\Bbb N}d(x(j),y(j))^p&=\sum_{l\in L_2}\sum_{j\in
K_2^l}d(r_j,s_j)^p\cr &<\sum_{l\in L_2}\sum_{j\in K_2^l}(\sqrt
2)^{-pl}d(u_j,v_j)^p\cr &<2C^p\sum_{l\in L_2}(\sqrt 2)^{-pl}\cr
&<+\infty,
\end{array}$$
so $(x,y)\in E((Z_{b_j})_{j\in\Bbb N},p)$. On the other hand, since
$2^{m_l}\ge N(l)$ we have
$$\begin{array}{ll}\sum_{j\in\Bbb N}\delta_q(T_j(x(j)),T_j(y(j)))^q&=\sum_{l\in L_2}\sum_{j\in
K_2^l}\delta_q(T_j(r_j),T_j(s_j))^q\cr &\ge\sum_{l\in L_2}\sum_{j\in
K_2^l}N(l)^{-q}\delta_q(T_j(u_j),T_j(v_j))^q\cr &>\sum_{l\in
L_2}\sum_{j\in K_2^l}N(l)^{-q}2^{qm_l}d(u_j,v_j)^p\cr &\ge
\sum_{l\in L_2}C^p\left(\frac{2^{m_l}}{N(l)}\right)^q\cr
&=+\infty,\end{array}$$ contradicting Lemma \ref{DH} again!
\hfill$\Box$

\begin{definition}
For two metric spaces $(X,d),(X',d')$ and $C,\alpha>0$. We say that
$X$ can {\it $C$-finitely H\"older$(\alpha)$ embed} into $X'$ if there
exists $A,D>0$ such that for every finite subset $F\subseteq X$,
there is $T_F:F\to X'$ satisfying, for $u,v\in F$,
\begin{enumerate}
\item[(1)] $d(u,v)\ge C\Rightarrow d'(T_F(u),T_F(v))\ge D$;
\item[(2)] $d(u,v)<C\Rightarrow A^{-1}d(u,v)^\alpha\le d'(T_F(u),T_F(v))\le
Ad(u,v)^\alpha$.
\end{enumerate}
While $\alpha=1$, we also say that $X$ can $C$-finitely Lipschitz embed
into $X'$.
\end{definition}

\begin{theorem} \label{freduce}
Let $(X,d)$ and $(Y_n,\delta_n),\,n\in\Bbb N$ be separable complete metric spaces,
$p,q\in[1,+\infty)$. If $X$ satisfies $({\rm link}(C))$ for some $C>0$, and
$F(X;p)\le_B E((Y_n)_{n\in\Bbb N};q)$, then
$X$ can $C$-finitely H\"older$(\frac{p}{q})$ embed into $\ell_q((Y_n)_{n\in\Bbb N},y^*)$ for any $y^*\in\prod_{n\in\Bbb N}Y_n$.
\end{theorem}

{\bf Proof.} Fix a sequence of finite subsets $F_n\subseteq
X,n\in\Bbb N$ such that
$$F_0\subseteq F_1\subseteq\cdots\subseteq F_n\subseteq\cdots$$
and $\bigcup_{n\in\Bbb N}F_n$ is dense in $X$. Let $(b_j)_{j\in\Bbb
N},(l_j)_{j\in\Bbb N}$ and $T_j:F_{b_j}\to \prod_{n=l_j}^{l_{j+1}-1}Y_n$ be from
the remarks before Lemma~\ref{D}. For convenience, we identify
$(\prod_{n=l_j}^{l_{j+1}-1}Y_n,\delta_q)$ with a subspace of $\ell_q((Y_n)_{n\in\Bbb N},y^*)$. Then
$T_j$ becomes a map $F_{b_j}\to\ell_q((Y_n)_{n\in\Bbb N},y^*)$.

Let us consider an arbitrary finite subset $F\subseteq X$. We can
find $k\in\Bbb N$ such that
\begin{enumerate}
\item[(a)] $k^{-1}\le d(u,v)$ for any distinct $u,v\in F$;
\item[(b)] $d(u,v)\le C-k^{-1}$ for any $u,v\in F$ with $d(u,v)<C$.
\end{enumerate}

For every $u\in F$, since $\bigcup_{j\in\Bbb N}F_{b_j}$ is dense in
$X$, there exists an $R(u)\in\bigcup_{j\in\Bbb N}F_{b_j}$ such that
$d(u,R(u))<(4k)^{-1}$. Then for any distinct $u,v\in F$, we have
$$d(R(u),R(v))<d(u,v)+(2k)^{-1}\le 2d(u,v),$$ and
$$d(R(u),R(v))>d(u,v)-(2k)^{-1}\ge\frac{1}{2}d(u,v).$$

From Lemmas~\ref{D} and \ref{fholder}, there exist $D>0,m\in\Bbb N$
and a sufficiently large $i$ such that $R(u)\in F_{b_i}$ for every
$u\in F$, and for $r,s\in F_{b_i}$,
\begin{enumerate}
\item[(i)] $d(r,s)\ge C-(2k)^{-1}\Rightarrow \delta_q(T_i(r),T_i(s))\ge D$;
\item[(ii)] $(2k)^{-1}\le d(r,s)<C\Rightarrow 2^{-m}d(r,s)^\frac{p}{q}\le \delta_q(T_i(r),T_i(s))\le
2^m d(r,s)^\frac{p}{q}$.
\end{enumerate}

We define $T_F:F\to\ell_q((Y_n)_{n\in\Bbb N},y^*)$ by $T_F(u)=T_i(R(u))$ for $u\in
F$.

For any $u,v\in F$ with $d(u,v)\ge C$, we have $d(R(u),R(v))\ge
C-(2k)^{-1}$. Then
$$\delta_q(T_F(u),T_F(v))=\delta_q(T_i(R(u)),T_i(R(v)))\ge D.$$

For any distinct $u,v\in F$ with $d(u,v)<C$, we have $k^{-1}\le
d(u,v)\le C-k^{-1}$. So $(2k)^{-1}\le d(R(u),R(v))\le
C-(2k)^{-1}<C$. Then
$$\begin{array}{ll}\delta_q(T_F(u),T_F(u))&=\delta_q(T_i(R(u)),T_i(R(v)))\cr
&\le 2^m d(R(u),R(v))^\frac{p}{q}\cr
&<2^{m+\frac{p}{q}}d(u,v)^\frac{p}{q},\end{array}$$ and
$$\begin{array}{ll}\delta_q(T_F(u),T_F(u))&=\delta_q(T_i(R(u)),T_i(R(v)))\cr
&\ge 2^{-m}d(R(u),R(v))^\frac{p}{q}\cr
&>2^{-(m+\frac{p}{q})}d(u,v)^\frac{p}{q}.\end{array}$$ Thus
$A=2^{m+\frac{p}{q}}$ and $D$ witness that $X$ can $C$-finitely
H\"older$(\frac{p}{q})$ embed into $\ell_q((Y_n)_{n\in\Bbb N},y^*)$. \hfill$\Box$

\begin{theorem}
Let $(X,d),(Y,\delta)$ be two separable complete metric spaces, $p,q\in[1,+\infty)$,
and let $Y_0\subseteq Y_1\subseteq Y_2\subseteq\cdots$ be a sequence of Borel subsets of $Y$ with $\bigcup_{n\in\Bbb N}Y_n$ dense in $Y$.
If $X$ satisfies $({\rm link}(C))$ for some $C>0$, then the following conditions are equivalent:
\begin{enumerate}
\item[(a)] $X$ can $C$-finitely H\"older$(\frac{p}{q})$ embed into $\ell_q((Y_n)_{n\in\Bbb N},y^*)$ for some $y^*\in\prod_{n\in\Bbb N}Y_n$.
\item[(b)] $F(X;p)\le_B E((Y_n)_{n\in\Bbb N};q)$.
\item[(c)] $F(X;p)\le_B F(Y;q)$.
\end{enumerate}
\end{theorem}

{\bf Proof.} Let $F_0\subseteq F_1\subseteq F_2\subseteq\cdots$ be a sequence of finite subsets of $X$ with $\bigcup_{n\in\Bbb N}F_n$ dense in $X$.

(a)$\Rightarrow$(b). Since $X$ can $C$-finitely
H\"older$(\frac{p}{q})$ embed into $\ell_q((Y_n)_{n\in\Bbb N},y^*)$, we can find
$A,D>0$, $T_n:F_n\to\ell_q((Y_n)_{n\in\Bbb N},y^*)$ such that, for $u,v\in F_n$,
\begin{enumerate}
\item[(1)] $d(u,v)\ge C\Rightarrow \delta_q(T_n(u),T_n(v))\ge D$;
\item[(2)] $d(u,v)<C\Rightarrow A^{-1}d(u,v)^\frac{p}{q}\le \delta_q(T_n(u),T_n(v))\le Ad(u,v)^\frac{p}{q}$.
\end{enumerate}
Then $F(X;p)\sim_B E((F_n)_{n\in\Bbb N};p)\le_B E((Y_n)_{n\in\Bbb N};q)$ follows from Theorem~\ref{reduce}.

(b)$\Rightarrow$(a) follows from Theorem \ref{freduce}.

(b)$\Rightarrow$(c). Let $(b_j)_{j\in\Bbb N},(l_j)_{j\in\Bbb N}$ and $T_j:F_{b_j}\to \prod_{n=l_j}^{l_{j+1}-1}Y_n$ be from
the remarks before Lemma~\ref{D}. Since every $F_{b_j}$ is finite, we can find finite subsets $U_n\subseteq Y_n$ for $l_j\le n<l_{j+1}$ such that
$T_j(u)\in\prod_{n=l_j}^{l_{j+1}-1}U_n$ for each $u\in F_{b_j}$. We can extend every $U_n$ to a finite subset $W_n\subseteq Y$ such that
$U_n\subseteq W_n$, $W_0\subseteq W_1\subseteq W_2\subseteq\cdots$ and $\bigcup_{n\in\Bbb N}W_n$ is dense in $Y$.

From Lemma \ref{D} with $C'=C$ and Lemma \ref{fholder} with $k=2^l$ , we can find $D>0$, $m\in\Bbb N$
and a strictly increasing sequence of natural numbers $(j_l)_{l\in\Bbb N}$ such that, for $r,s\in F'_l\stackrel{\rm Def}{=}F_{b_{j_l}}$, we have
\begin{enumerate}
\item[(i)] $d(r,s)\ge C\Rightarrow\delta_q(T_{j_l}(r),T_{j_l}(s))\ge D$;
\item[(ii)] $2^{-l}\le d(r,s)<C\Rightarrow 2^{-m}d(r,s)^{\frac{p}{q}}\le \delta_q(T_{j_l}(r),T_{j_l}(s))\le 2^m d(r,s)^{\frac{p}{q}}$.
\end{enumerate}

Then Corollary \ref{reduce1} gives
$$F(X;p)\sim_B E((F'_l)_{l\in\Bbb N};p)\le_B E((W_n)_{n\in\Bbb N};q)\sim_B F(Y;q).$$

(c)$\Rightarrow$(b). Find a sequence of finite subsets $V_n\subseteq Y_n,\,n\in\Bbb N$ such that $V_0\subseteq V_1\subseteq V_2\subseteq\cdots$ and
$\bigcup_{n\in\Bbb N}V_n$ is dense in $Y$. Then we have $F(X;p)\le_B F(Y;q)\sim_B E((V_n)_{n\in\Bbb N};q)\le_B E((Y_n)_{n\in\Bbb N};q)$.
\hfill$\Box$

\begin{corollary}
Let $X,Y$ be two separable complete metric spaces, $p,q\in[1,+\infty)$. If $X$ satisfies $({\rm link}(C))$ for some $C>0$, then
the following conditions are equivalent:
\begin{enumerate}
\item[(a)] $X$ can $C$-finitely H\"older$(\frac{p}{q})$ embed into $\ell_q(Y,y^*)$ for some $y^*\in Y^\Bbb N$.
\item[(b)] $F(X;p)\le_B E(Y;q)$.
\item[(c)] $F(X;p)\le_B F(Y;q)$.
\end{enumerate}
\end{corollary}

\end{document}